\documentclass[12pt]{amsart}
\usepackage{amscd, amssymb}

\newcommand{\bdes}{\begin{description}}
\newcommand{\edes}{\end{description}}
\newcommand{\beqn}{\begin{eqnarray*}}
\newcommand{\eeqn}{\end{eqnarray*}}
\newcommand{\PP }{{\mathbb P}}
\newcommand{\QQ }{{\mathbb Q}}
\newcommand{\CC }{{\mathbb C}}
\newcommand{\ZZ }{{\mathbb Z}}

\newcommand{\hb }{{\hbar }}
\newcommand{\ti }{\times}

\newcommand{\ra }{\rightarrow}

\begin{document}

\author{Bumsig Kim}
\title{Quantum Hyperplane Section Principle For Concavex
Decomposable Vector Bundles}
\maketitle
\pagestyle{plain}

\section{Introduction}

The Lefschetz hyperplane section theorem says that
there is an intimate relationship between the
cohomology group of an ambient space and that of
a smooth zero locus of a positive
line bundle over the ambient space.
Roughly speaking, it's quantum version says that there is an
intimate  relationship between quantum
cohomology rings of the ambient space and that of 
a smooth zero locus of the decomposable spanned vector bundle
\cite{BCKV1, BCKV2}.

In paper \cite{Kim}, we proved the quantum analog 
when the 
ambient space is a generalized flag manifold and
the decomposable vector bundle is convex. In it we claimed
that the quantum analog can be generalized to the case
when the bundle is concavex and decomposable.
We explain the claim in this paper. 
As an application we reprove the multiple cover formula.

This work is originally motivated
by the Lian-Liu-Yau paper \cite{LLY}. A mirror theorem
for concave bundle spaces over symplectic toric manifolds
is established by Givental \cite{GL}.
There is also a work
of Elezi's paper \cite{Elezi} in the generalization of Givental's work
for concavex decomposable vector bundle over projective spaces.

The result of the paper was announced in Bay Area
Joint Symplectic Geometry
Seminar in April, 1998 at Stanford University.

\section{Notation}

\begin{itemize}
\item Let $X$ be a generalized flag manifold.

\item Let
$V$ be a vector bundle decomposable to the 
direct sum of convex or concave line bundles $L_j$, $j=1,...,k$.
A line
bundle $L$ is called convex if 
$H^1(\PP ^1, f^*(L_j))=0$ 
for any morphism $f: \PP ^1\ra X$.
If $H^0(\PP ^1, f^*(L_j))=0$ for any morphism $f: \PP ^1\ra X$,
$L$ is called concave.  
We call $V$ concavex and decomposable following \cite{LLY}.

\item Let $p_i$,  $i=1,...,l$, be the divisor classes of $X$ such that 
$\bigoplus _{i=1}^l\ZZ_{\ge 0} p_i$ is the closed K\"ahler cone.

\item Let $q_i$, $i=1,...,l$, denote indeterminants and
$q^{\beta} := \prod _{i=1}^l q_i^{<p_i,\beta >}$, where
$\beta$ is an effective class of $H_2(X,\ZZ )$.
So, $q^\beta \in \QQ [q_1,...,q_l]$.

\item Let $T':=\CC ^\ti$, the complex torus, which
acts on $V$ fiberwise only using the scalar product of
the vector spaces of fibers. Let $T'$ act on the base space $X$
trivially. So $V$ is a $T'$-equivariant vector bundle. 

\item Denote by $H^*_{(T')}$ the rational field of the equivariant cohomology
     $H^*_{T'}(point, \QQ )$ of one point.  Let $H^*_{(T')}(X):=
    H^*(X,\QQ )\otimes _{\QQ}H^*_{(T')}$

\item Given an effective class $\beta \in H_2(X,\ZZ )$,
$$\overline{M}_{0,n}(X,\beta )$$ denotes the moduli space of equivalence
classes $$[(C,f;x_1,...,x_n)]$$ of $n$-marked
stable maps of genus 0 and type $\beta$. An element of it can be represented
by a holomorphic map $f$ to $X$  from a connected nodal curve $C$
of arithmetic genus zero with
$f_*([C])=\beta\in H_2(X)$ and $n$-distinct ordered
 nonsingular points $x_i\in C$.
The moduli space consists of stable ones. 
The space is a compact complex orbifold with dimension
$\dim X + <c_1(X),\beta>+n-3$ since $T_X$ is generated by
global sections.


\end{itemize}

\section{Quantum Cohomology associated with $(X,V)$}

Associated to the $T'$-equivariant bundle $V=\bigoplus _{j=1}^k L_j$ 
over $X$, we define a new quantum
product on $H^*_{(T')}(X)$ following \cite{GE}. 
Modified Poincar\' e pairings based on $V$ will be utilized.

Let $A$ and $B$ denote equivariant classes in $H^*_{(T')}(X)$.  
Define $<,>_{0}^V$, a nondegenerate inner product in $ H^*_{(T')}(X)$
by

$$ <A,B>_0^V:=\int _X A B E_{T'}(V),$$ where 
$E_{T'}(V):=\prod _i E_{T'}(L_i)$ and
\[ E_{T'}(L) := \left\{
    \begin{array}{ll} Euler _{T'}(L) & \mbox{if $L$ is convex} \\
                     Euler _{T'}(L)^{-1} & \mbox{if $L$ is concave}.
\end{array} \right. \]

Here we use the $T'$-equivariant Euler classes $Euler _{T'}(L)$, 
so that the class $E_{T'}(V)$
is invertible over the coefficient ring $H^*_{(T')}$ and 
thus $<,>^V_0$ is indeed nondegenerate.

Introduce
the induced $T'$-equivariant vector (orbi-)bundles \[
[L]_\beta := \left\{
    \begin{array}{ll}R^0\pi _*(ev_{N+1})^* (L)  & \mbox{if $L$ is convex} \\
                   R^1\pi _*(ev_{N+1})^* (L)   & \mbox{if $L$ is concave}
\end{array} \right. , \] where
$ev_{N+1}$ denotes the evaluation map at ($N+1$)-th marked points
from $\overline{M}_{0,N+1}(X,\beta )$ to $X$ and $\pi$ denotes
the forgetting-last-marked-point map from $\overline{M}_{0,N+1}(X,\beta )$
to $\overline{M}_{0,N}(X,\beta )$. 

Let $A_i\in H^*_{(T')}(X)$, $i=1,...,N$.
Define 
$<,...,>_{\beta}^V$, $N$-correlators, by

$$ <A_1,...,A_N>_\beta ^V=\int _{\overline{M}_{0,N}(X,\beta )}
   ev_1^*(A_1)... ev_N^*(A_N) Euler _{T'} (V_\beta ) $$
where 
      \[ V_\beta := \bigoplus _i [L_i]_\beta \]
    In this
definition if $\beta =0$, then we assume $N\ge 3$.
Here $T'$ acts $V_\beta$ fiberwise and
\[ ev_i: \overline{M}_{0,N}(X,\beta ) \ra X, \ \ i=1,...,N \] 
are the evaluation maps at $i$-th marked points.

With these N- correlators and the nondegenerate pairing $<,>^L_0$,
one can define a big/small quantum cohomology on $H^*_{(T')}(X)$ and also
a quantum differential system, its fundamental solution, and so on.
For instance, let $A$ and $B$ be in $H^*_{(T')}(X)$, then
the small quantum product
$$A*_{V}B\in H^*_{(T')}(X)\otimes _{\QQ}\QQ [[q]]$$
is defined by the requirement 
\[ <A*_V B, C>_0^V = \sum _{\beta} q^{\beta} <A,B,C>_{\beta}^V, \]
for all $C\in H^*_{(T')}(X)$. So, $H^*_{(T')}(X)\otimes _{\QQ}\QQ [[q_1,...,
q_l]]$ has the small quantum ring structure based on $V$ and $V_\beta$.

Introduce formal parameters $t_1,...,t_l$ and the relations $q_i=e^{t_i}$.
The quantum differential system is a formal family of formal first order
partial differential equations in $t_1,...,t_l$ 
with $\hb$ as a formal parameter:
\[ p_i*_V f(t,q) = \hb\frac {\partial}{\partial t_i}f(t,q) \]
where $f(t,q)\in H ^*_{(T')}(X)[\hb ^{-1}][t_1,...,t_l][[q_1,...,q_l]]$.
Here we treat $q_i\frac{\partial }{\partial q_i}
=\frac{\partial}{\partial t_i}$ 
formally.

{\em Let
$(ev_N)_{*,V}$ be defined as the adjoint of the pullback $ev_N^*$ with
respect to the {\bf new}
Poincar\' e pairings on the N-marked moduli spaces and $X$.} So, by
the very definition of the pushforward,
$$\int _{\overline{M}_{0,N}(X,\beta )} A\cup ev_N^*(B)\cup Euler _{T'}(V_\beta) 
= \int _X (ev_N)_{*,V}(A)\cup B\cup E_{T'}(V)$$
for $A\in H_{(T')}^*(\overline{M}_{0,N}(X,\beta ))$ and $B\in H_{(T')}^*(X)$.

Now we describe a fundamental solution to the quantum differential system.
For any given $A\in H^*_{T'}(X)$,
\[ f_A(t,q) := \sum _{\beta \ne 0} q^{\beta}
(ev_2)_{*,V}(ev_1^*(A)\frac{\exp ({ev_1^*(pt)/\hb})}{\hb -c} )
+A\exp (pt/\hb) \]
is a solution where $pt:=\sum _i p_it_i$
and $c$ is the nonequivariant first Chern class of
the universal cotangent line bundle at first marked points.
To show it one may use WDVV, string and divisor equations \cite{GL}
or use divisor equation and topological recursion relation \cite{Pa}.

\section{The Givental correlator}

Now define the so-called {\bf Givental correlator}
\[ J_{\beta}^V:= ev_{*, V}(\frac{1}{\hbar (\hbar - c)}), \]
where $ev: \overline{M}_{0,1}(X,\beta ) \ra X $ is the evaluation
map, and let
$J^V_0:=1$.

It is obtained from special components from solutions
to the small quantum differential equations. That is, 
$$<J^V_\beta , \exp (pt/\hb )A>^V_0
=\int _{\overline{M}_{0,2}(X,\beta )} 
\frac{\exp (ev_1^*(pt)/\hb) ev_1^*(A)}{\hb -c} \cup ev_2^*(1) 
\cup Euler _{T'}(V_\beta)$$
if $\beta\ne 0$. 

\bigskip

On the other hand,
for $\beta \neq 0$,

\[ <J^V_\beta , 1>_0^V =
\frac{-2}{\hb ^3}\int _{\overline{M}_{0,0}(X,\beta)}Euler_{T'}(V_\beta)
+o(\hb ^{-3}),\] where the integral over the no-marked moduli space
of $Euler(V_\beta)$ is sometimes very interesting. Some examples
are as follows.

\bigskip

{\bf Example 1}
If $X=\PP ^4$ and $V=\mathcal{O}(5)$, then
the integration computes the degree of the virtual fundamental
class of degree $\beta$ for a smooth quintic. In turn it provides
the integer numbers of \lq\lq almost" rational curves of given homotopic
types in the quintic defined by Ruan.

{\bf Example 2} If $X=\PP ^1$ and $V={\mathcal O}(-1)\oplus
    {\mathcal O}(-1)$, the integral is the multiple covering
contribution.

\bigskip

When $V$ is the rank zero bundle (that is, there are no Euler things
in correlators)
we denote
\[ J_{\beta}^X:= ev _* (\frac{1}{\hbar (\hbar - c)}). \] 

\bigskip

{\bf Example 3}
When $X$ is a projective space $\PP ^n$,
  \[ J^{\PP ^n}
 = 1+\sum _{d=1}^{\infty}
q^d\frac 1{((p+\hb )... (p+ d\hb ))^{n+1}} \] as in \cite{GE}.

\bigskip

\section{Quantum Hyperplane Section Principle}
{\bf We want to compare $J_{\beta}^V$ and $J_{\beta}^X$.} 

Define
 \[ H^L_{\beta}:= \left\{
    \begin{array}{lll} 
         \prod _{m=1}^{<c_1(L),\beta >}    (c^{T'}_1(L)+m\hb )  
    & \mbox{if $L$ is convex} \\ & \\
 \prod ^{0}_{m=<c_1(L),\beta >+1}  (c^{T'}_1(L)+m\hb ) 
   & \mbox{if $L$ is concave}.
    \end{array} \right. \]
Let $$ H^V_\beta = \prod _i H^{L_i}_\beta $$ if $\beta\ne 0$
and $H^V_0=1$.
{\bf It will be called the correcting Euler class for $V$.}
Here $c^{T'}_1(L)$ is the $T'$-equivariant first Chern class of $L$.

Define
$$J^V(q_1,...,q_l)
:= \sum _{\beta\in H_2(X,\ZZ)} q^\beta J_\beta ^V $$
and $$I^V (q_1,...,q_l):=\sum _{\beta \in H_2(X,\ZZ)}
q^\beta J_\beta ^X H_\beta ^V .$$ 

The degree of $q_i$ is uniquely defined by the requirement: 
$$c_1(TX)-\sum _{\text{convex }L_i}c_1(L_i)+\sum _{\text {concave } L_i}
c_1(L_i)
= \sum _i(\deg q_i) p_i .$$

\bigskip

{\bf Theorem} {\em Suppose each $\deg q_i$ is nonnegative.
Then 
  $$J^V=
e^{f_0+f_{-1}/\hb+\sum p_i f_i/\hb}
I^V (q_1e^{f_1},...,q_le^{f_l})$$
for unique $q$-series $f_i$ without constant terms 
where $\deg f_i= 0$ for $i=0,...,l$
and $\deg f_{-1}=1$.
In particular, 
if $I^V=1+O(\hb ^{-2})$, then $J^V=I^V$.}

\bigskip 
{\bf Proof}. The proof in \cite{Kim} for the convex case 
works for this general, concavex case, word for word.

\bigskip

{\bf Example of multiple cover formula}
In this case it is easy to see that $I^V$ starts with $1+O(\hb ^{-2})$
in the expansion in $\hb ^{-1}$.
Therefore, $J^V=I^V$ and thus $<J^V_d,1>^V_0=<I^V_d,1>^V_0
=\int _{\PP ^1}\frac{1}{(p+d\hb)^2}=-2\frac{1}{\hb ^3 d^3}$.
We obtain the multiple cover formula which is first proven by Manin.

\vspace{+10 pt}
\noindent
Department of Mathematics \\
University of California - Davis \\
Davis, CA 95616\\
bumsig@math.ucdavis.edu
\end{document}